# ERGODIC PROPERTIES OF POISSONIAN ID PROCESSES

### By Emmanuel Roy

### *Université Paris 13*

We show that a stationary IDp process (i.e., an infinitely divisible stationary process without Gaussian part) can be written as the independent sum of four stationary IDp processes, each of them belonging to a different class characterized by its Lévy measure. The ergodic properties of each class are, respectively, nonergodicity, weak mixing, mixing of all order and Bernoullicity. To obtain these results, we use the representation of an IDp process as an integral with respect to a Poisson measure, which, more generally, has led us to study basic ergodic properties of these objects.

**1. Introduction.** A stochastic process is said to be *infinitely divisible* (ID) if, for any positive integer $k$, it equals, in distribution, the sum of $k$ independent and identically distributed processes. These processes are fundamental objects in probability theory, the most popular being the intensively studied Lévy processes (see, e.g., [19]). We will focus here on ID stationary processes $\{X_n\}_{n \in \mathbb{Z}}$. Stationary Gaussian processes have a particular place among stationary ID processes and have already been the subject of very deep studies (see [7] for recent results). We will concentrate on non-Gaussian ID processes; Maruyama [8] first started their study. Since the late eighties many authors are looking for criteria of ergodicity, weak mixing or mixing of a general ID process, exhibiting examples, studying particular sub-families [mainly symmetric $\alpha$-stable ($S\alpha S$) processes]. We mention the result of Rosiński and Żak [17] which shows the equivalence of ergodicity and weak mixing for general ID processes. Some *factorizations* (for the convolution product) have been obtained in the $S\alpha S$ case, in particular, Rosiński [13] has shown that a $S\alpha S$ process can be written in a unique way as the independent sum of three $S\alpha S$ processes, one being called *mixed moving* *average* (which is mixing), the second *harmonizable* (nonergodic) and the third not in the aforementioned categories and which is potentially the most interesting (see [15]) (note that Rosiński has developed, in [14], a multidimensional version of this factorization). Recently, this third part has been split by Pipiras and Taqqu (see [12]) and Samorodnitsky managed to isolate (through a factorization) the "maximal" ergodic component of a $S\alpha S$ process (see [18]). Factorizations already appeared in [9], where the ID objects were ID point processes.

The fundamental tool in the study of a non-Gaussian ID process is its *Lévy measure*. In the stationary case, its existence has been shown by Maruyama in [8]: it is a (shift-)stationary measure on $\mathbb{R}^{\mathbb{Z}}$, which might be infinite, related to the distribution of the ID process by the characteristic functions of its finite-dimensional distributions through an extended *Lévy–Khintchine* formula. A general ID process is the independent sum of a Gaussian process and a *Poissonian* (IDp) process, the latter being uniquely determined by its Lévy measure. Reciprocally, if we are given a (shift-)stationary measure on $\mathbb{R}^{\mathbb{Z}}$, under some mild conditions, it can be seen as the Lévy measure of a unique IDp stationary process.

Our main result consists in establishing the following factorization result: every IDp stationary process can be written in a unique way as the independent sum of four IDp processes which are, respectively, nonergodic, weakly mixing, mixing (of all order) and Bernoulli (Theorem 5.5 and Proposition 5.7).

The proof is divided in several steps which have their own interest. The first step is based on the following remark: if the support of the Lévy measure can be partitioned into invariant sets, then the restrictions to these sets of the measure are the Lévy measures of processes that form a factorization of the initial process. We point out here that it may happen that a stationary ID process can be factorizable into infinitely many components, however, we only consider factorizations that make sense in terms of ergodic behavior of each class. It is remarkable that those distinct behaviors are naturally linked to those of the corresponding Lévy measures. Thus, it is essential to get a better understanding of general dynamical systems (particularly with infinite measure) and to study decompositions along their invariant sets. Section 2 presents some elements of ergodic theory. In particular, we recall a decomposition, mostly due to Hopf, Krengel and Sucheston (see [6]), of an invariant measure into the sum four invariant measures which are the restrictions of the initial measure to as many invariant sets with distinctive properties (Proposition 2.11). Section 3 presents some basic facts of spectral theory that will be used later. There are no new results in Sections 2 and 3.

Then, back to Lévy measures, we have to link the different categories to the corresponding ergodic properties of the underlying ID process. To do so, we use the representation due to Maruyama [8] of any IDp process as a stochastic integral with respect to the Poisson measure with the Lévy measure as intensity. In ergodic terms, we will say that an IDP process is a *factor* of the Poisson suspension constructed above its Lévy measure. We thus are led to a specific study of Poisson suspensions built above dynamical systems that is the subject of Section 4. This study is mostly based upon the particular structure of the associated $L^2$-space, which admits a *chaotic* decomposition: the Fock factorization of the $L^2$-space associated to the underlying dynamical system. This preliminary work allows us to elucidate absence of ergodicity, weak mixing and mixing of all order of a Poisson suspension. We also give a criterion for the Bernoulli property.

In Section 5 we first recall the basic facts on infinitely divisible processes and then apply the results of the preceding sections to their Lévy measure. Thanks to our factorization, ergodic properties can be easily derived. In Section 6 we give an explicit form of all stationary IDp processes with a dissipative Lévy measure. In cases where the process is square integrable, some spectral criteria for ergodic behaviors can be established (Section 7).

In Section 8 were the $\alpha$-stable case is treated, we show that our factorization preserves the distributional properties, that is, each of the four components is $\alpha$-stable. We can thus replace in this context the previously obtained factorization of Rosiński [13], as well as the refinements of Pipiras and Taqqu [12] and Samorodnitsky [18].

**2. Elements of ergodic theory.** Let $(\Omega, \mathcal{F}, \mu)$ be a $\sigma$-finite Lebesgue space in the following sense: there exists a probability measure $\nu$, equivalent to $\mu$, such that $(\Omega, \mathcal{F}, \nu)$ is a Lebesgue space in its traditional acceptation. Let $T$ be a bijective bimeasurable transformation that preserves $\mu$. The quadruplet $(\Omega, \mathcal{F}, \mu, T)$ is called *dynamical system*, or *system* for short.

The aim of this section is to introduce basic notions and terminology used in the study of dynamical systems. We first concentrate on the structure of a general dynamical system that will lead us to the decomposition in Proposition 2.11 which is a compilation of known results. The rest of the section is devoted to notions specific to dynamical systems with a probability measure. The book of Aaronson [1] covers most of the definitions and results exposed here.

In the following, if $\phi$ is a measurable map defined on $(\Omega, \mathcal{F}, \mu, T)$, the image measure of $\mu$ by $\phi$ is denoted $\phi^\star(\mu)$.

2.1. *Factors, isomorphic systems.* Consider another dynamical system $(\Omega', \mathcal{F}', \mu', T')$.

DEFINITION 2.1. Call $(\Omega', \mathcal{F}', \mu', T')$ a *factor* of $(\Omega, \mathcal{F}, \mu, T)$ if there exists a map $\varphi$, measurable from $(\Omega, \mathcal{F})$ to $(\Omega', \mathcal{F}')$ such that $\varphi^\star(\mu) = \mu'$ and $\varphi \circ T = T' \circ \varphi$. If $\varphi$ is invertible, then $(\Omega, \mathcal{F}, \mu, T)$ and $(\Omega', \mathcal{F}', \mu', T')$ are said to be *isomorphic*.

2.2. *Ergodicity.*

DEFINITION 2.2.    The *invariant σ-field* of $(\Omega, \mathcal{F}, \mu, T)$ is the sub-σ-field $\mathcal{I}$ of $\mathcal{F}$ that contains the sets $A \in \mathcal{F}$ such that $T^{-1}A = A$ ($A$ is said to be *T-invariant*).

This definition leads to the following one:

DEFINITION 2.3.    $(\Omega, \mathcal{F}, \mu, T)$ is said to be *ergodic* if, for all set $A \in \mathcal{I}$,

$$\mu(A) = 0 \quad \text{or} \quad \mu(A^c) = 0.$$

2.3. *Dissipative and conservative transformations.*

DEFINITION 2.4.    A set $A \in \mathcal{F}$ is called a *wandering set* if the sets $\{T^{-n}A\}_{n \in \mathbb{Z}}$ are disjoint.

We denote by $\mathfrak{D}$ the (measurable) union of all the wandering sets for $T$, this set is $T$-invariant. Its complement is denoted by $\mathfrak{C}$.

DEFINITION 2.5.    We call $(\Omega, \mathcal{F}, \mu, T)$ *dissipative* if $\mathfrak{D} = \Omega$ mod. $\mu$. If $\mathfrak{C} = \Omega$ mod. $\mu$, then $(\Omega, \mathcal{F}, \mu, T)$ is said *conservative*.

LEMMA 2.6.    *There exists a wandering set $W$ such that $\mathfrak{D} = \bigcup_{n \in \mathbb{Z}} T^{-n}W$ mod. $\mu$.*

PROPOSITION 2.7 (Hopf decomposition).    *The Hopf decomposition is the partition $\{\mathfrak{D}, \mathfrak{C}\}$.*
    $(\Omega, \mathcal{F}, \mu_{|\mathfrak{D}}, T)$ *is dissipative and* $(\Omega, \mathcal{F}, \mu_{|\mathfrak{C}}, T)$ *is conservative.*

2.4. *Type* $\mathbf{II}_1$ *and type* $\mathbf{II}_\infty$.    The following proposition is a consequence of the decomposition found in [1], page 47.

PROPOSITION 2.8.    *Let $(\Omega, \mathcal{F}, \mu, T)$ be a dynamical system. There exists a unique partition $\{\mathfrak{P}, \mathcal{N}\}$ of $\Omega$ in $T$-invariant sets such that there exists a $T$-invariant probability measure equivalent to $\mu_{|\mathfrak{P}}$ and that there does not exist a nonzero $T$-invariant probability measure absolutely continuous with respect to $\mu_{|\mathcal{N}}$. We have $\mathfrak{P} \subset \mathfrak{C}$ and $\mathfrak{D} \subset \mathcal{N}$. $(\Omega, \mathcal{F}, \mu_{|\mathfrak{P}}, T)$ is said to be of type* $\mathbf{II}_1$ *and* $(\Omega, \mathcal{F}, \mu_{|\mathcal{N}}, T)$ *of type* $\mathbf{II}_\infty$.

REMARK.    We use the notion of type $\mathbf{II}_\infty$ in an abusive manner since it includes dissipative transformations. However it is very convenient in our context.

2.5. *Zero type and positive type.*

**Definition 2.9.** Let $(\Omega, \mathcal{F}, \mu, T)$ be a dynamical system.

$(\Omega, \mathcal{F}, \mu, T)$ is said to be of *zero type* if, for all $A \in \mathcal{F}$ such that $0 < \mu(A) < +\infty$, $\mu(A \cap T^{-k}A) \to 0$ as $k$ tends to $+\infty$.

$(\Omega, \mathcal{F}, \mu, T)$ is said to be of *positive type* if, for all $A \in \mathcal{F}$ such that $\mu(A) > 0$, $\overline{\lim}_{k \to \infty} \mu(A \cap T^{-k}A) > 0$.

**Remark.** By using similar arguments as in Theorem 5.5, page 58 in [11], it is easy to see that $(\Omega, \mathcal{F}, \mu, T)$ is of zero type if and only if, for all $A, B \in \mathcal{F}$ such that $0 < \mu(A) < +\infty$ and $0 < \mu(B) < +\infty$, $\mu(A \cap T^{-k}B) \to 0$ as $k$ tends to $+\infty$.

Krengel and Sucheston obtained the following decomposition (see [6], page 155):

**Proposition 2.10.** *There exists a partition $\{\mathcal{N}_0, \mathcal{N}_+\}$ of $\Omega$ in $T$-invariant sets such that $(\Omega, \mathcal{F}, \mu_{|\mathcal{N}_0}, T)$ [resp. $(\Omega, \mathcal{F}, \mu_{|\mathcal{N}_+}, T)$] is of zero type (resp. of positive type). We have $\mathfrak{D} \subset \mathcal{N}_0$ and $\mathfrak{P} \subset \mathcal{N}_+ \subset \mathfrak{C}$.*

Note that Aaronson in [1] calls positive part, the part of type $\mathbf{II}_1$ and null part, the part of type $\mathbf{II}_\infty$.

We can group all these decompositions in the following proposition:

**Proposition 2.11** (Canonical decomposition). *Let $(\Omega, \mathcal{F}, \mu, T)$ be a dynamical system. By defining $\mu_B := \mu_{|\mathfrak{D}}$, $\mu_m := \mu_{|\mathcal{N}_0 \cap \mathfrak{C}}$, $\mu_{wm} := \mu_{|\mathcal{N}_+ \cap \mathcal{N}}$ and $\mu_{ne} := \mu_{|\mathfrak{P}}$ (this choice of notation is motivated by Theorem 4.8), we can write, in a unique way,*

$$\mu = \mu_B + \mu_m + \mu_{wm} + \mu_{ne},$$

*where:*

$(\Omega, \mathcal{F}, \mu_B, T)$ *is dissipative.*

$(\Omega, \mathcal{F}, \mu_m, T)$ *is conservative of zero type.*

$(\Omega, \mathcal{F}, \mu_{wm}, T)$ *is of positive and $\mathbf{II}_\infty$ type.*

$(\Omega, \mathcal{F}, \mu_{ne}, T)$ *is of type $\mathbf{II}_1$.*

**Remark.** Note that none of these categories is empty, [5] provides various examples of conservative type $\mathbf{II}_\infty$ dynamical systems.

2.6. *The case of a probability measure.*   We assume here that $\mu(\Omega) = 1$.

THEOREM 2.12 (Birkhoff ergodic theorem).   *Let $f \in L^1(\mu)$, then, $\mu$-a.e. and in $L^1(\mu)$*

$$\lim_{n \to \infty} \frac{1}{n} \sum_{k=1}^{n} f \circ T^k = \mu(f|\mathcal{I}),$$

*where $\mu(f|\mathcal{I})$ is the conditional expectation of $f$ with respect to the invariant $\sigma$-algebra.*

DEFINITION 2.13.   $(\Omega, \mathcal{F}, \mu, T)$ is said to be *weakly mixing* if, for all $A, B \in \mathcal{F}$,

$$\lim_{n \to \infty} \frac{1}{n} \sum_{k=1}^{n} |\mu(A \cap T^{-k}B) - \mu(A)\mu(B)| = 0.$$

$(\Omega, \mathcal{F}, \mu, T)$ is said to be *mixing of order $m$* if, for all $A_1, \ldots, A_m \in \mathcal{F}$, for all strictly increasing sequences of integer $n_{1,k}, \ldots, n_{m,k}$,

$$\lim_{k \to \infty} |\mu(T^{n_{1,k}} A_1 \cap \cdots \cap T^{n_{1,k} + \cdots + n_{m,k}} A_m) - \mu(A_1) \cdots \mu(A_m)| = 0.$$

$(\Omega, \mathcal{F}, \mu, T)$ is said to be *mixing* if it is mixing of order 2, that is, if for all $A, B \in \mathcal{F}$,

$$\lim_{n \to \infty} |\mu(A \cap T^{-n}B) - \mu(A)\mu(B)| = 0.$$

We now introduce a dynamical system that will constantly be used in the paper. We consider here the space $\mathbb{R}^{\mathbb{Z}}$ of $\mathbb{Z}$-indexed sequences. The natural $\sigma$-algebra is the product $\sigma$-algebra $\mathcal{B}^{\otimes \mathbb{Z}}$, where $\mathcal{B}$ is the natural Borel $\sigma$-algebra on $\mathbb{R}$. The transformation is the shift $T$ that acts in the following way:

$$T\{x_i\}_{i \in \mathbb{Z}} = \{x_{i+1}\}_{i \in \mathbb{Z}}.$$

The dynamical system $(\mathbb{R}^{\mathbb{Z}}, \mathcal{B}^{\otimes \mathbb{Z}}, \mu, T)$ is the canonical space of the stationary process of distribution $\mu$.

DEFINITION 2.14.   The system associated to an i.i.d. process is called a *Bernoulli scheme*. A system $(\Omega, \mathcal{F}, \mu, T)$ is said to be *Bernoulli* if it is isomorphic to a Bernoulli scheme.

We end this section by the following proposition:

PROPOSITION 2.15.   *We have the implications:*

*Bernoulli $\Rightarrow$ mixing of order $n \Rightarrow$ mixing $\Rightarrow$ weakly mixing $\Rightarrow$ ergodic.*

*Moreover, these six properties are shared by all the factors.*

**3. Spectral theory.** Here we only give results that will be needed in the rest of the paper. See [2] and [1] for details and proofs.

3.1. *Hilbert space, unitary operator and spectral measure.* We consider a complex Hilbert space $(H, \langle \cdot \rangle)$ endowed with a unitary operator $U$. To each vector $f \in H$, we can associate a finite measure $\sigma_f$ on $[-\pi, \pi[$, called the *spectral measure of $f$* by the formula

$$\hat{\sigma}_f(n) := \langle U^n f, f \rangle = \int_{[-\pi,\pi[} e^{inx} \sigma_f(dx).$$

Let $C(f)$ be the closure of the linear space generated by the family $\{U^n f\}_{n \in \mathbb{Z}}$, $C(f)$ is called the *cyclic space* of $f$. We summarize the following properties in the following proposition:

PROPOSITION 3.1. *There exists an isomorphism $\phi$ between $C(f)$ and $L^2(\sigma_f)$ with $\phi(f) = 1$ and such that the unitary operator $h \mapsto e^{i \cdot} h$ on $L^2(\sigma_f)$ is conjugate to $U$ by $\phi$.*

3.2. *Maximal spectral type.* On $(H, \langle \cdot \rangle, U)$ there exists a finite measure $\sigma_M$ such that, for all $f \in H$, $\sigma_f \ll \sigma_M$. The (equivalence class of the) measure $\sigma_M$ is called the *maximal spectral type* of $(H, \langle \cdot \rangle, U)$. Moreover, for all finite measures $\sigma \ll \sigma_M$, there exists a vector $g$ such that $\sigma_g = \sigma$.

3.3. *Application to ergodic theory.* A dynamical system $(\Omega, \mathcal{F}, \mu, T)$ induces a complex Hilbert space, the space $L^2(\mu)$ endowed with a unitary operator $U : f \mapsto f \circ T$.

3.3.1. *The case of a probability measure.* We restrict the study to the orthocomplement of the constant functions in $L^2(\mu)$. That is, we note $L_0^2(\mu) := L^2(\mu) \ominus \mathbb{C}\langle 1 \rangle$ and we call *reduced maximal spectral type* of $(\Omega, \mathcal{F}, \mu, T)$ the maximal spectral type of $(L_0^2(\mu), U)$. We recover the following ergodic properties on the reduced maximal spectral type $\sigma_M$:

PROPOSITION 3.2. *$(\Omega, \mathcal{F}, \mu, T)$ is ergodic if and only if $\sigma_M\{0\} = 0$.*
*$(\Omega, \mathcal{F}, \mu, T)$ is weakly mixing if and only if $\sigma_M$ is continuous.*
*$(\Omega, \mathcal{F}, \mu, T)$ is mixing if and only if $\sigma_M$ is a Rajchman measure [i.e., $\hat{\sigma}_f(n) \to 0$ as $|n|$ tends to $+\infty$].*

3.3.2. *The infinite measure case.* Since constant nonzero functions are not in $L^2(\mu)$, we do not impose the restriction made in the preceding section. $\mathbf{II}_\infty$ and zero types are spectral properties:

PROPOSITION 3.3. $(\Omega, \mathcal{F}, \mu, T)$ *is of type* $\mathbf{II}_\infty$ *if and only if* $\sigma_M$ *is continuous and this condition is also equivalent to* $\sigma_M\{0\} = 0$.
$(\Omega, \mathcal{F}, \mu, T)$ *is of zero type if and only if* $\sigma_M$ *is Rajchman.*

PROOF. The fact that $(\Omega, \mathcal{F}, \mu, T)$ is of type $\mathbf{II}_\infty$ if and only if $\sigma_M$ is continuous can be found in [1], page 74.

We now prove that $\sigma_M\{0\} = 0$ implies that $\sigma_M$ is continuous. Assume that $\sigma_M$ is not continuous, then $(\Omega, \mathcal{F}, \mu, T)$ is not of type $\mathbf{II}_\infty$, that is, there exists a $T$-invariant probability measure $\nu$ such that $\nu \ll \mu$. The function $\sqrt{\frac{d\nu}{d\mu}}$ is in $L^2(\mu)$ and, since it is $T$-invariant, its spectral measure is the Dirac mass at 0.

The proof of the last statement on zero type systems is completely similar to the mixing case for probability preserving systems (see, e.g., pages 57–58 in [11]). □

## 4. Poisson suspensions.

In this section we will recall basic facts on the intensively studied Poisson measures, which are random discrete measures on an underlying measure space. The particular case we are interested in, that is, when the distribution of the Poisson measure is preserved by a well chosen transformation (and then called Poisson suspension), has received much less attention ([2] provides a few pages on Poisson suspensions and references, mainly under the scope of statistical mechanics). The particular form, in chaos, of the $L^2$-space associated to the Poisson suspension allows a useful spectral analysis similar to the Gaussian case.

4.1. *Definitions.* We consider a $\sigma$-finite Lebesgue space $(\Omega, \mathcal{F}, \mu)$. Let $\{A_n\}_{n \in \mathbb{N}}$ be a countable measurable partition of $\Omega$ such that $\mu(A_n) < \infty$ for all $n \in \mathbb{N}$ and let $(M_\Omega, \mathcal{M}_\mathcal{F})$ be the space of measures $\nu$ on $(\Omega, \mathcal{F})$ satisfying $\nu(A_n) \in \mathbb{N}$ for all $n \in \mathbb{N}$. $\mathcal{M}_\mathcal{F}$ is the smallest $\sigma$-algebra on $M_\Omega$ such that the mappings $\nu \to \nu(A)$ are measurable for all $A \in \mathcal{F}$ of finite $\mu$-measure. We denote by $N$ the identity on $(M_\Omega, \mathcal{M}_\mathcal{F})$.

DEFINITION 4.1. We call *Poisson measure* the triplet $(M_\Omega, \mathcal{M}_\mathcal{F}, \mathcal{P}_\mu)$, where $\mathcal{P}_\mu$ is the unique probability measure such that, for all finite collections $\{A_i\}$ of elements belonging to $\mathcal{F}$, disjoint and of finite $\mu$-measure, the $\{N(A_i)\}$ are independent and distributed as the Poisson law of parameter $\mu(A_i)$. The underlying space $(\Omega, \mathcal{F}, \mu)$ will be called the *base*.

Assume now that $T$ is an invertible and measure preserving transformation on $(\Omega, \mathcal{F}, \mu)$; it is easily verified that the map $T^\star$ defined on $M_\Omega$ by $T^\star(\nu) = \nu \circ T^{-1}$ is also a bijective transformation which preserves the probability $\mathcal{P}_\mu$.

DEFINITION 4.2. The dynamical system $(M_\Omega, \mathcal{M}_\mathcal{F}, \mathcal{P}_\mu, T^\star)$ is called the *Poisson suspension* above the base $(\Omega, \mathcal{F}, \mu, T)$.

4.2. *Product structure.* The independence properties of a Poisson suspension along invariant subsets imply the following:

LEMMA 4.3. *Let $(\Omega, \mathcal{F}, \mu, T)$ be a dynamical system and suppose there exists a partition $\{\Omega_i\}_{1 \leq i \leq k}$ of $\Omega$ into $k$ $T$-invariant sets of nonzero $\mu$-measure.*

*Then $(M_\Omega, \mathcal{M}_\mathcal{F}, \mathcal{P}_\mu, T^\star)$ is isomorphic to the direct product*

$$(M_\Omega^k, \mathcal{M}_\mathcal{F}^{\otimes k}, \mathcal{P}_{\mu_{|\Omega_0}} \otimes \cdots \otimes \mathcal{P}_{\mu_{|\Omega_k}}, T^\star \times \cdots \times T^\star).$$

4.3. *General $L^2$ properties of a Poisson suspension.* In this section we recall the basic facts on the Fock space structure of the $L^2$-space associated to a Poisson measure $(M_\Omega, \mathcal{M}_\mathcal{F}, \mathcal{P}_\mu)$. Section 10.4 in [10] is a reference for this section.

4.3.1. *Fock factorization.*

DEFINITION 4.4. The *Fock factorization* of the Hilbert space $K$ is the Hilbert space $\mathbf{exp}K$ given by

$$\mathbf{exp}K := \mathfrak{S}^0 K \oplus \mathfrak{S}^1 K \oplus \cdots \oplus \mathfrak{S}^n K \oplus \cdots,$$

where $\mathfrak{S}^n K$ is the $n$th symmetric tensor power of $K$ and is called the *$n$th chaos*, with $\mathfrak{S}^0 K = \mathbb{C}$.

On $\mathbf{exp}K$, the set of *exponential vectors* is particularly important

$$\mathcal{E}_h := 1 \oplus h \oplus \frac{1}{\sqrt{2!}}(h \otimes h) \oplus \cdots \oplus \frac{1}{\sqrt{n!}}(h \otimes \cdots \otimes h) \oplus \cdots$$

for $h \in K$.

They form a linearly dense part in $\mathbf{exp}K$ and satisfy the identity

$$\langle \mathcal{E}_h, \mathcal{E}_g \rangle_{\mathbf{exp}K} = \exp \langle h, g \rangle_K.$$

Now suppose we are given an operator $U$ on $K$ with norm at most 1, it extends naturally to an operator $\tilde{U}$ on $\mathbf{exp}K$ called the *exponential* of $U$, by acting on each chaos via the formula

$$\tilde{U}(h \otimes \cdots \otimes h) = (Uh) \otimes \cdots \otimes (Uh)$$

leading to the identity,

$$\tilde{U}\mathcal{E}_h = \mathcal{E}_{Uh}.$$

4.3.2. *Fock space structure of $L^2(\mathcal{P}_\mu)$.* Call $\Delta_n$ the diagonal in $\Omega^n$ (the $n$-uplets with identical coordinates). Multiple integrals, for $f$ in $L^1(\mu) \cap L^2(\mu)$, are defined by

$$J^{(n)}(f)$$
$$:= \int \cdots \int_{\Delta_n^c} f(x_1) \cdots f(x_n)(N(dx_1) - \mu(dx_1)) \cdots (N(dx_n) - \mu(dx_n)).$$

THEOREM 4.5.   *There exists an isometry between $L^2(\mathcal{P}_\mu)$ and $\exp[L^2(\mu)]$ mapping $J^{(n)}(f)$ to $\sqrt{n!} \underbrace{f \otimes \cdots \otimes f}_{n\ times}$ for any $n \geq 1$ and $f$ in $L^1(\mu) \cap L^2(\mu)$.*

We thus have the isometry formula:

$$\langle J^{(n)}(f), J^{(p)}(g) \rangle_{L^2(\mathcal{P}_\mu)} = n!(\langle f, g \rangle_{L^2(\mu)})^n 1_{n=p}.$$

Call $\mathfrak{H}$ the set of functions $h$, finite linear combination of indicator functions of elements of $\mathcal{F}$ with finite $\mu$-measure. Through the natural isometry, the exponential vectors $\mathcal{E}_h$ are

$$\mathcal{E}_h(\nu) = \exp\left(-\int_\Omega h \, d\mu\right) \prod_{x \in \nu}(1 + h(x)).$$

They form a linearly dense part in $L^2(\mathcal{P}_\mu)$, moreover, $\mathbb{E}_{\mathcal{P}_\mu}[\mathcal{E}_h] = 1$.

4.4. *Spectral properties of a Poisson suspension.* We now consider the case of a dynamical system $(\Omega, \mathcal{F}, \mu, T)$ and its associated Poisson suspension $(M_\Omega, \mathcal{M}_\mathcal{F}, \mathcal{P}_\mu, T^\star)$. It is obvious that the unitary operator $f \mapsto f \circ T^\star$ acting on $L^2(\mathcal{P}_\mu)$ is the exponential of the corresponding unitary operator on $L^2(\mu)$, $g \mapsto g \circ T$. From this simple remark, it can be deduced, as in the Gaussian case (see Chapter 14 in [2] for details), with the following very important properties:

PROPOSITION 4.6.   *If $\sigma_M$ is the maximal spectral type of $(\Omega, \mathcal{F}, \mu, T)$, then on the $n$th chaos, the maximal spectral type of $U$ is $\sigma_M^{*n}$. The (reduced) maximal spectral type of the Poisson suspension $(M_\Omega, \mathcal{M}_\mathcal{F}, \mathcal{P}_\mu, T^\star)$ is $e(\sigma_M) := \sum_{n \geq 1} \frac{1}{n!} \sigma_M^{*n}$.*

4.5. *Ergodic properties of a Poisson suspension.* In this section we consider a system $(M_\Omega, \mathcal{M}_\mathcal{F}, \mathcal{P}_\mu, T^\star)$, where $\mu = \mu_B + \mu_m + \mu_{wm} + \mu_{ne}$ from the decomposition in Proposition 2.11. Lemma 4.3 immediately implies the following:

PROPOSITION 4.7.   *$(M_\Omega, \mathcal{M}_\mathcal{F}, \mathcal{P}_\mu, T^\star)$ is isomorphic to*

$$(M_\Omega^4, \mathcal{M}_\mathcal{F}^{\otimes 4}, \mathcal{P}_{\mu_B} \otimes \mathcal{P}_{\mu_m} \otimes \mathcal{P}_{\mu_{wm}} \otimes \mathcal{P}_{\mu_{ne}}, T^\star \times T^\star \times T^\star \times T^\star).$$

We now look at the ergodic properties in each class:

THEOREM 4.8. $(M_\Omega, \mathcal{M}_\mathcal{F}, \mathcal{P}_{\mu_{ne}}, T^\star)$ *is not ergodic.*
$(M_\Omega, \mathcal{M}_\mathcal{F}, \mathcal{P}_{\mu_{wm}}, T^\star)$ *is weakly mixing, not mixing.*
$(M_\Omega, \mathcal{M}_\mathcal{F}, \mathcal{P}_{\mu_m}, T^\star)$ *is mixing of all orders.*
$(M_\Omega, \mathcal{M}_\mathcal{F}, \mathcal{P}_{\mu_B}, T^\star)$ *is Bernoulli.*

PROOF. Since $(\Omega, \mathcal{F}, \mu_{ne}, T)$ is not of type $\mathbf{II}_\infty$, from Proposition 3.3, its maximal spectral type has an atom at 0 and this implies that, thanks to Proposition 4.6, this atom at 0 is part of the (reduced) maximal spectral type of $(M_\Omega, \mathcal{M}_\mathcal{F}, \mathcal{P}_{\mu_{ne}}, T^\star)$ and thus prevents ergodicity.

The fact that $(M_\Omega, \mathcal{M}_\mathcal{F}, \mathcal{P}_{\mu_{wm}}, T^\star)$ is weakly mixing is a direct consequence of the successive application of Propositions 3.3, 4.6 and 3.2. Since $\sigma_M$ is not Rajchman, it cannot be mixing.

If now we consider $(\Omega, \mathcal{F}, \mu_m, T)$, this system is of zero type, that is to say, for all $A \in \mathcal{F}$, $B \in \mathcal{F}$ of finite $\mu$-measure, $\mu_m(A \cap T^{-k}B)$ tends to 0 as $k$ tends to infinity.

We are going to generalize the identity $\langle \mathcal{E}_h, \mathcal{E}_g \rangle_{L^2(\mathcal{P}_{\mu_m})} = \exp\langle h, g \rangle_{L^2(\mu_m)}$:

$$\mathbb{E}_{\mathcal{P}_{\mu_m}}[\mathcal{E}_{h_1}\mathcal{E}_{h_2}\cdots\mathcal{E}_{h_n}]$$
$$= \exp \sum_{1 \le i_1 < i_2 \le n} \int h_{i_1} h_{i_2} \, d\mu_m + \cdots$$
$$+ \sum_{1 \le i_1 < i_2 < \cdots < i_n \le n} \int h_{i_1} \cdots h_{i_n} \, d\mu_m.$$

We show, more generally, the following formula for functions $h_1, \ldots, h_n$ of $\mathfrak{H}$:

$$\mathcal{E}_{h_1}\mathcal{E}_{h_2}\cdots\mathcal{E}_{h_n}$$
$$= \mathcal{E}_{(1+h_1)(1+h_2)\cdots(1+h_n)-1} \exp \sum_{1 \le i_1 < i_2 \le n} \int h_{i_1} h_{i_2} \, d\mu_m + \cdots$$
$$+ \sum_{1 \le i_1 < i_2 < \cdots < i_n \le n} \int h_{i_1} \cdots h_{i_n} d\mu_m.$$

At rank 2, the computation is easy; let $n \ge 2$ and suppose that the formula is true at this rank.

Let $h_1, \ldots, h_n, h_{n+1}$ be functions in $\mathfrak{H}$.

We first evaluate $\mathcal{E}_{(1+h_1)(1+h_2)\cdots(1+h_n)-1}\mathcal{E}_{h_{n+1}}$. The formula, at rank 2, gives us

$$\mathcal{E}_{(1+h_1)(1+h_2)\cdots(1+h_n)-1}\mathcal{E}_{h_{n+1}}$$

$$= \exp \int h_{n+1}((1+h_1)(1+h_2)\cdots(1+h_n)-1)\,d\mu_m$$

$$\times \mathcal{E}_{(1+h_1)(1+h_2)\cdots(1+h_n)(1+h_{n+1})-1}.$$

But $\exp \int h_{n+1}((1+h_1)(1+h_2)\cdots(1+h_n)-1)\,d\mu_m$ equals

$$\exp \sum_{i=1}^{n} \int h_i h_{n+1}\,d\mu_m + \cdots$$

$$+ \sum_{1 \le i_1 < i_2 < \cdots < i_n \le n} \int h_{i_1}\cdots h_{i_n} h_{n+1}\,d\mu_m.$$

Combining this result with the formula at rank $n$, we show that the formula is true at rank $n+1$ and this ends the proof by recurrence.

To show mixing of order $n$ with the functions $\mathcal{E}_{h_1},\ldots,\mathcal{E}_{h_n}$ with $h_1,\ldots,h_n$ in $\mathfrak{H}$, take $n$ strictly increasing sequences of integers $p_{1,k},\ldots,p_{n,k}$ and denote by $a_{i,k} := p_{1,k} + \cdots + p_{i,k}$. We have to show that
$\mathbb{E}_{\mathcal{P}_{\mu_m}}[\mathcal{E}_{h_1}\circ T^{\star a_{1,k}}\mathcal{E}_{h_2}\circ T^{\star a_{2,k}}\cdots \mathcal{E}_{h_n}\circ T^{\star a_{n,k}}]$ tends to

$$\mathbb{E}_{\mathcal{P}_{\mu_m}}[\mathcal{E}_{h_1}]\cdots \mathbb{E}_{\mathcal{P}_{\mu_m}}[\mathcal{E}_{h_n}] = 1.$$

But

$$\mathbb{E}_{\mathcal{P}_{\mu_m}}[\mathcal{E}_{h_1}\circ T^{\star a_{1,k}}\mathcal{E}_{h_2}\circ T^{\star a_{2,k}}\cdots \mathcal{E}_{h_n}\circ T^{\star a_{n,k}}]$$

$$= \mathbb{E}_{\mathcal{P}_{\mu_m}}[\mathcal{E}_{h_1\circ T^{a_{1,k}}}\mathcal{E}_{h_2\circ T^{a_{2,k}}}\cdots \mathcal{E}_{h_n\circ T^{a_{n,k}}}]$$

and then, from the preceding formula, we have to show that quantities of the kind $\int h_i \circ T^{a_{i,k}}\cdots h_j \circ T^{a_{j,k}}\,d\mu_m$, $i < j$, tend to 0.

The functions $h_i$ are finite linear combinations of indicator functions of sets of finite $\mu$-measure, then, expanding the integral $\int h_i \circ T^{a_{i,k}}\cdots h_j \circ T^{a_{j,k}}\,d\mu_m$, we obtain a finite linear combination of quantities of the kind $\mu_m(T^{-a_{i,k}}A_l \cap \cdots \cap T^{-a_{j,k}}A_m)$. But these quantities tend to 0 since

$$\mu_m(T^{-a_{i,k}}A_l \cap \cdots \cap T^{-a_{j,k}}A_m) \le \mu_m(T^{-a_{i,k}}A_l \cap T^{-a_{j,k}}A_m)$$

and

$$\mu_m(T^{-a_{i,k}}A_l \cap T^{-a_{j,k}}A_m) = \mu_m(A_l \cap T^{p_{i+1,k}+\cdots+p_{j,k}}A_m).$$

We thus have the mixing of order $n$ on the exponential vectors $\mathcal{E}_{h_1},\ldots,\mathcal{E}_{h_n}$, and, by standard approximation arguments, taking advantage of the properties of these vectors, we get mixing of order $n$ for the suspension.

$(\Omega,\mathcal{F},\mu_B,T)$ is dissipative, so, from Lemma 2.6, there exists a wandering set $W$ such that $\Omega = \bigcup_{n\in\mathbb{Z}} T^{-n}W$ mod. $\mu_B$. Denote by $\mathcal{W}$ the $\sigma$-field generated by $A \in \mathcal{F}$ such that $A \subset W$. Then $\mathcal{M}_{\mathcal{W}}$ generates $\mathcal{M}_{\mathcal{F}}$ (i.e., $\mathcal{M}_{\mathcal{F}} = \bigvee_{n\in\mathbb{Z}} T^{\star-n}\mathcal{M}_{\mathcal{W}}$) and, thanks to the independence properties of a Poisson measure, the $\sigma$-fields $T^{\star-n}\mathcal{M}_{\mathcal{W}}$ are independent. Hence, $(M_{\Omega},\mathcal{M}_{\mathcal{F}},\mathcal{P}_{\mu_B},T^{\star})$ is Bernoulli.  $\square$

REMARK. The content of this theorem is apparently due to Marchat in his Ph.D. dissertation as pointed out by Grabinski in [4], we have heard of Grabinski's paper, which is cited nowhere, at the "Galley proofs" stage of the preparation of this document.

A direct consequence of this theorem is that a Poisson suspension is ergodic (and weakly mixing) if and only if the base is $\mathbf{II}_\infty$. This has also been proved in [3] which contains also results of modern ergodic theory on Poisson suspensions.

**5. Infinitely divisible stationary processes.** After a few generalities on stationary processes, we next introduce the notion of infinite divisibility for these processes which is an immediate generalization of the finite-dimensional case (the book of K. Sato [19] is a reference on this vast subject). The accompanying tools such as the Lévy measure find its equivalent notion for processes as shown by Maruyama in [8]. This measure is the key object that will allow us to connect results of the preceding sections to prove Theorem 5.5, which was the motivation for this work, and to deduce their ergodic properties in Theorem 5.7.

5.1. *Dynamical system associated to a stationary stochastic process.* We consider $(\mathbb{R}^{\mathbb{Z}}, \mathcal{B}^{\otimes \mathbb{Z}}, \mu, T)$ introduced in Section 2.6, $\mu$ may be infinite. When we will deal with stationary processes, only the measure will change throughout the study and, to simplify, we will often use it to designate such a system. Affirmations such as "$\mu$ is ergodic" or "$\mu$ is dissipative" will be shortening of "$(\mathbb{R}^{\mathbb{Z}}, \mathcal{B}^{\otimes \mathbb{Z}}, \mu, T)$ is ergodic" or "$(\mathbb{R}^{\mathbb{Z}}, \mathcal{B}^{\otimes \mathbb{Z}}, \mu, T)$ is dissipative." We will try to keep the notation $X := \{X_0 \circ T^n\}_{n \in \mathbb{Z}}$ for the identity on $(\mathbb{R}^{\mathbb{Z}}, \mathcal{B}^{\otimes \mathbb{Z}})$, $X_0$ being the "coordinate at 0" map $\{x_i\}_{i \in \mathbb{Z}} \mapsto x_0$. $X$, $\{X_n\}_{n \in \mathbb{Z}}$, $\{X_0 \circ T^n\}_{n \in \mathbb{Z}}$, $\mu$ or $(\mathbb{R}^{\mathbb{Z}}, \mathcal{B}^{\otimes \mathbb{Z}}, \mu, T)$ is essentially the same object.

5.2. *Convolution of processes.* We consider the mapping "sum" with values in $(\mathbb{R}^{\mathbb{Z}}, \mathcal{B}^{\otimes \mathbb{Z}})$ which associates $\{x_i + y_i\}_{i \in \mathbb{Z}}$ to $(\{x_i\}_{i \in \mathbb{Z}}, \{y_i\}_{i \in \mathbb{Z}})$. Given two distributions $\mathbb{P}_1$ and $\mathbb{P}_2$ on $(\mathbb{R}^{\mathbb{Z}}, \mathcal{B}^{\otimes \mathbb{Z}})$, we call $\mathbb{P}_1 * \mathbb{P}_2$ the "convolution of $\mathbb{P}_1$ with $\mathbb{P}_2$." $\mathbb{P}_1 * \mathbb{P}_2$ is the image distribution of $\mathbb{P}_1 \otimes \mathbb{P}_2$ by the mapping already defined. Since this operation is clearly associative, we can denote $\mathbb{P}^{*k}$ to be the convolution of $k$ identical copies of $\mathbb{P}$.

DEFINITION 5.1. Let $\mathbb{P}$ be a distribution on $(\mathbb{R}^{\mathbb{Z}}, \mathcal{B}^{\otimes \mathbb{Z}})$, $\mathbb{P}$ is *infinitely divisible* (ID) if, for all integers $k$, there exists a distribution $\mathbb{P}_k$ on $(\mathbb{R}^{\mathbb{Z}}, \mathcal{B}^{\otimes \mathbb{Z}})$ such that $\mathbb{P} = \mathbb{P}_k^{*k}$.

We remark that this definition forces the finite-dimensional distributions to be ID.

5.2.1. *Lévy measure of an ID stationary process.* We have, as in the finite-dimensional case, a representation, due to Maruyama (see [8]), of characteristic functions of the finite-dimensional distributions of an ID stationary process of distribution $\mathbb{P}$ (we denote by $a$ a sequence $\{a_i\}_{i \in \mathbb{Z}}$ where only a finite number of coordinates are nonzero and call $\mathfrak{A}$ their union in $\mathbb{R}^{\mathbb{Z}}$):

$$(5.1) \quad \begin{aligned} &\mathbb{E}[\exp i \langle a, X \rangle] \\ &\quad = \exp\left[-\tfrac{1}{2}\langle Ra, a \rangle + i\langle a, b_\infty \rangle + \int_{\mathbb{R}^{\mathbb{Z}}} (e^{i\langle a,x \rangle} - 1 - i\langle c(x), a \rangle) Q(dx)\right], \end{aligned}$$

where $R$ is the covariance function of a centered stationary Gaussian process, $b_\infty \in \mathbb{R}^{\mathbb{Z}}$ is a sequence identically equal to $b$ and $Q$ is a $\sigma$-finite measure on $(\mathbb{R}^{\mathbb{Z}}, \mathcal{B}^{\otimes \mathbb{Z}})$ invariant with respect to the shift and such that $Q\{0\} = 0$ (where $\{0\}$ is the identically zero sequence), $\int_{\mathbb{R}^{\mathbb{Z}}} (x_0^2 \wedge 1) Q(dx) < +\infty$ and $c(x)_i = -1_{]-\infty,-1[} + x_i 1_{[-1,1]} + 1_{]1,\infty[}$.

$\langle R, b, Q \rangle$ is called the *generating triplet* of $\mathbb{P}$.

The dynamical system $(\mathbb{R}^{\mathbb{Z}}, \mathcal{B}^{\otimes \mathbb{Z}}, Q, T)$ will be our main concern in the sequel.

When the process is integrable and centered, we have the following representation, where $R$ and $Q$ are unchanged:

$$(5.2) \quad \mathbb{E}[\exp i \langle a, X \rangle] = \exp\left[-\tfrac{1}{2}\langle Ra, a \rangle + \int_{\mathbb{R}^{\mathbb{Z}}} (e^{i\langle a,x \rangle} - 1 - i\langle a, x \rangle) Q(dx)\right].$$

Finally, if the process only takes positive values (and then without Gaussian part), we can write down its finite-dimensional distribution through their Laplace transforms, with $a \in \mathfrak{A} \cap \mathbb{R}_+^{\mathbb{Z}}$:

$$(5.3) \quad \mathbb{E}[\exp -\langle a, X \rangle] = \exp\left[-\langle a, b_\infty \rangle - \int_{\mathbb{R}^{\mathbb{Z}}} (1 - e^{-\langle a,x \rangle}) Q(dx)\right].$$

If, moreover, it is integrable, under this representation, we have

$$\mathbb{E}[X_0] = b + \int_{\mathbb{R}^{\mathbb{Z}}} x_0 Q(dx).$$

REMARK 5.2. If we are given a covariance function $R$, a drift $b$ and a measure $Q$ satisfying the hypothesis specified above, it determines the distribution of an ID process of generating triplet $\langle R, b, Q \rangle$ by defining its finite-dimensional distribution through the representation (5.1). Then we can apprehend the extraordinary variety of the process at our disposal.

DEFINITION 5.3. An ID process is said to be *Poissonian* (IDp) if its generating triplet does not possess a Gaussian part.

In the sequel, when we will speak of IDp process with Lévy measure $Q$, we will consider a process whose generating triplet is $\langle 0, 0, Q \rangle$ under the representation (5.1). Of course, the drift has no impact in our study.

### 5.3. *First examples and representation.*

5.3.1. *Canonical example.* Maruyama in [8] has given the canonical example of an IDp stationary process:

We consider a Poisson suspension $(M_\Omega, \mathcal{M}_\mathcal{F}, \mathcal{P}_\mu, T^\star)$ above $(\Omega, \mathcal{F}, \mu, T)$ and a real function $f$ defined on $(\Omega, \mathcal{F}, \mu, T)$ such that $\int_\Omega \frac{f^2}{1+f^2} \, d\mu < +\infty$. We define the stochastic integral $I(f)$ by the limit in probability, as $n$ tends toward infinity, of

$$\int_{|f|>1/n} f \, dN - \int_{|f|>1/n} c(f) \, d\mu.$$

Then the process $X = \{I(f) \circ T^{\star n}\}_{n \in \mathbb{Z}}$ is IDp and its distribution is given by

$$\mathbb{E}[\exp i\langle a, X\rangle] = \exp\left[\int_\Omega \exp\left(i\sum_{n\in\mathbb{Z}} a_n f \circ T^n\right) - 1 - i\sum_{n\in\mathbb{Z}} a_n c(f \circ T^n) \, d\mu\right]$$

for $a \in \mathfrak{A}$.

Maruyama has also shown in [8] that all the IDp processes can be represented this way: consider $Q$, the Lévy measure of an IDp process of generating triplet $\langle 0, 0, Q\rangle$, let $(M_{\mathbb{R}^\mathbb{Z}}, \mathcal{M}_{\mathcal{B}^{\otimes\mathbb{Z}}}, \mathcal{P}_Q, T^\star)$ be the Poisson suspension with base $(\mathbb{R}^\mathbb{Z}, \mathcal{B}^{\otimes\mathbb{Z}}, Q, T)$ and $f$ the mapping $X_0: \{x_i\}_{i\in\mathbb{Z}} \mapsto x_0$.

**Theorem 5.4** (Maruyama). *The process $\{I(X_0) \circ T^{\star n}\}_{n\in\mathbb{Z}}$ admits $\langle 0, 0, Q\rangle$ as generating triplet.*

This theorem is crucial since it allows us to consider an IDp process as a factor of a Poisson suspension, precisely the Poisson suspension constructed above its Lévy measure. The factor map is $\nu \to \{I(x_0) \circ T^{\star n}\}_{n\in\mathbb{Z}}$, defined on $M_{\mathbb{R}^\mathbb{Z}}$ with values in $\mathbb{R}^\mathbb{Z}$.

### 5.4. *First factorization.*

It is obvious that the convolution of two ID distributions is still ID, the class of this type of distributions being closed under convolution. Given a stationary ID distribution, we ask when it is *factorizable*, that is, can it be written as the convolution of two or more ID distributions? An immediate factorization comes from the representation (5.1):

Suppose that $\mathbb{P}$ admits the triplet $\langle R, b, Q\rangle$. If $\mathbb{P}_s$ admits the triplet $\langle sR, sb, sQ\rangle$ and $\mathbb{P}_{1-s}$ admits the triplet $\langle (1-s)R, (1-s)b, (1-s)Q\rangle$ with $0 < s < 1$, then $\mathbb{P} = \mathbb{P}_s * \mathbb{P}_{1-s}$.

The representation (5.1) allows another more interesting factorization. Letting $\mathbb{P}_R$ of triplet $\langle R, 0, 0\rangle$, $\mathbb{P}_b$ of triplet $\langle 0, b, 0\rangle$ and $\mathbb{P}_Q$ of triplet $\langle 0, 0, Q\rangle$, we have

$$\mathbb{P} = \mathbb{P}_R * \mathbb{P}_b * \mathbb{P}_Q,$$

where $\mathbb{P}_R$ is the distribution of a stationary centered Gaussian process, $\mathbb{P}_b$ is the distribution of a constant process and $\mathbb{P}_Q$ the distribution of an IDp process.

5.5. *Factorization through invariant components of the Lévy measure.* We can apply to $Q$ the decomposition $Q = Q_B + Q_m + Q_{wm} + Q_{ne}$ along the four disjoint shift-invariant subsets as in Proposition 2.11. By considering (5.1), we get the following factorization result:

THEOREM 5.5 (Factorization of a stationary IDp process).  *Let $\mathbb{P}$ be the distribution of a stationary IDp process. $\mathbb{P}$ can be written in the unique way:*

$$\mathbb{P} = \mathbb{P}_{Q_B} * \mathbb{P}_{Q_m} * \mathbb{P}_{Q_{wm}} * \mathbb{P}_{Q_{ne}},$$

*where:*
  $(\mathbb{R}^{\mathbb{Z}}, \mathcal{B}^{\otimes \mathbb{Z}}, Q_B, T)$ *is dissipative,*
  $(\mathbb{R}^{\mathbb{Z}}, \mathcal{B}^{\otimes \mathbb{Z}}, Q_m, T)$ *is conservative of zero type,*
  $(\mathbb{R}^{\mathbb{Z}}, \mathcal{B}^{\otimes \mathbb{Z}}, Q_{wm}, T)$ *is of type $\mathbf{II}_\infty$ and of positive type,*
  $(\mathbb{R}^{\mathbb{Z}}, \mathcal{B}^{\otimes \mathbb{Z}}, Q_{ne}, T)$ *is of type $\mathbf{II}_1$.*

Since these classes are not empty for the corresponding Poisson suspensions, we deduce they are not empty for the IDp processes by considering stochastic integrals with respect to these Poisson suspensions.

5.6. *Ergodic properties of stationary IDp processes.* Before enunciating the properties of each class, we will need the following lemma which is the interpretation, in our framework, of a computation done by Rosiński and Żak in [17]. Their computation led to show that, if $X$ is an IDp process, the spectral measure of $e^{iX_0} - \mathbb{E}[e^{iX_0}]$ has the form $|\mathbb{E}[e^{iX_0}]|^2 e(m)$ (we still use the notation

$$e(m) := \sum_{k=1}^{+\infty} \frac{1}{k!} m^{*k},$$

where $m$ is a finite measure on $[-\pi, \pi[$. We will see that $m$ is indeed itself a spectral measure, but for the system associated to the Lévy measure of $X$.

LEMMA 5.6.  *Let $X$ be an IDp process of Lévy measure $Q$. The spectral measure of $e^{iX_0} - \mathbb{E}[e^{iX_0}]$ is $|\mathbb{E}[e^{iX_0}]|^2 e(\sigma)$, where $\sigma$ is the spectral measure of $e^{iX_0} - 1$ under $Q$.*

PROOF.  In [17], the following formula is established:

$$\mathbb{E}[e^{iX_0} \overline{e^{iX_k}}] = |\mathbb{E}[e^{iX_0}]|^2 \left( \exp\left[ \int_{\mathbb{R}^2} (e^{ix} - 1)(\overline{e^{iy} - 1}) Q_{0,k}(dx, dy) \right] \right),$$

where $Q_{0,k}$ is the Lévy measure of the ID vector $(X_0, X_k)$. But, since we make use of the Lévy measure of processes, this formula can be written into

$$\mathbb{E}[e^{iX_0} \overline{e^{iX_k}}] = |\mathbb{E}[e^{iX_0}]|^2 \left( \exp\left[ \int_{\mathbb{R}^2} (e^{ix_0} - 1)(\overline{e^{ix_k} - 1}) Q(dx) \right] \right),$$

which equals

$$|\mathbb{E}[e^{iX_0}]|^2 (\exp \hat{\sigma}(k)) = |\mathbb{E}[e^{iX_0}]|^2 \left( \sum_{n=0}^{+\infty} \frac{1}{n!} (\hat{\sigma}(k))^n \right),$$

where $\sigma$ is the spectral measure of $e^{iX_0} - 1$ under $Q$. The conclusion follows. $\square$

THEOREM 5.7.   $(\mathbb{R}^{\mathbb{Z}}, \mathcal{B}^{\otimes \mathbb{Z}}, \mathbb{P}_{Q_{ne}}, T)$ *is not ergodic.*
$(\mathbb{R}^{\mathbb{Z}}, \mathcal{B}^{\otimes \mathbb{Z}}, \mathbb{P}_{Q_{wm}}, T)$ *is weakly mixing.*
$(\mathbb{R}^{\mathbb{Z}}, \mathcal{B}^{\otimes \mathbb{Z}}, \mathbb{P}_{Q_m}, T)$ *is mixing of all order.*
$(\mathbb{R}^{\mathbb{Z}}, \mathcal{B}^{\otimes \mathbb{Z}}, \mathbb{P}_{Q_B}, T)$ *has the Bernoulli property.*

PROOF.   There exists a probability measure $\nu$ which is $T$-invariant and equivalent to $Q_{ne}$. Let $f := \sqrt{\frac{dQ_{ne}}{d\nu}}$ [note that $\frac{dQ_{ne}}{d\nu}$ is just $(\frac{d\nu}{dQ_{ne}})^{-1}$] and $\lambda \in \mathbb{R}$.

The spectral measure of $e^{i\lambda X_0} - 1$ under $Q_{ne}$ is the spectral measure of $fe^{i\lambda X_0} - f$ under $\nu$. The set $\{f < a\}$ is $T$-invariant since $f$ is $T$-invariant, moreover, this set is of nonzero measure if $a$ is large enough. Thus, the spectral measure of $fe^{i\lambda X_0} - f$ under $\nu$ is the sum of the spectral measures of $(fe^{i\lambda X_0} - f)1_{\{f < a\}}$ and $(fe^{i\lambda X_0} - f)1_{\{f \geq a\}}$ under $\nu$.

If $(fe^{i\lambda X_0} - f)1_{\{f < a\}}$ is centered, we have

$$\int_{\mathbb{R}^{\mathbb{Z}} \cap \{f < a\}} f(x) e^{i\lambda x_0} \nu(dx) = \int_{\mathbb{R}^{\mathbb{Z}} \cap \{f < a\}} f(x) \nu(dx) \in \mathbb{R}.$$

This implies

$$\int_{\mathbb{R}^{\mathbb{Z}} \cap \{f < a\}} f(x)[1 - \cos(\lambda x_0)] \nu(dx) = 0.$$

Since $f$ is nonnegative on $\{f < a\}$ $\nu$-a.e., this implies that $\cos(\lambda X_0) = 1$ on $\{f < a\}$ $\nu$-a.e. or that $\lambda X_0 = 0$ mod. $\pi$. But this is impossible for all $\lambda \in \mathbb{R}$ simultaneously.

That is, there exists $\lambda \in \mathbb{R}$ such that $(fe^{i\lambda X_0} - f)1_{\{f < a\}}$ is not centered and this implies that the spectral measure of $e^{i\lambda X_0} - 1$ under $Q_{ne}$ possesses an atom at 0. This atom is also in the spectral measure of $e^{i\lambda X_0} - \mathbb{E}[e^{i\lambda X_0}]$ by Lemma 5.6 and then in the maximal spectral type, which prevents ergodicity.

$(\mathbb{R}^{\mathbb{Z}}, \mathcal{B}^{\otimes \mathbb{Z}}, \mathbb{P}_{Q_{wm}}, T)$ is a factor of $(M_{\mathbb{R}^{\mathbb{Z}}}, \mathcal{M}_{\mathcal{B}^{\otimes \mathbb{Z}}}, \mathcal{P}_{Q_{wm}}, T^*)$, which is weakly mixing.

The rest of the properties are proved in the same way by considering the system as a factor of the corresponding Poisson suspension whose properties, such as mixing of all order and Bernoullicity, are inherited by its factors (see Proposition 2.15).   □

We are now able to give a new proof of the important theorem of Rosiński and Żak (see [17]).

THEOREM 5.8.   *If $\mathbb{P}$ is IDp and ergodic, then $\mathbb{P}$ is weakly mixing.*

PROOF.   Let $\mathbb{P} = \mathbb{P}_{Q_B} * \mathbb{P}_{Q_m} * \mathbb{P}_{Q_{wm}} * \mathbb{P}_{Q_{ne}}$ be the factorization of $\mathbb{P}$ from Theorem 5.5 with $\mathbb{P}_{Q_{ne}}$ nontrivial. Thus, $\{X_n\}_{n \in \mathbb{Z}}$ of distribution $\mathbb{P}$ can be seen as the independent sum of $\{X_n^1\}_{n \in \mathbb{Z}}$ of distribution $\mathbb{P}_{Q_B} * \mathbb{P}_{Q_m} * \mathbb{P}_{Q_{mm}} * \mathbb{P}_{Q_{wm}}$ and $\{X_n^2\}_{n \in \mathbb{Z}}$ of distribution $\mathbb{P}_{Q_{ne}}$. From the first part of the proof of Theorem 5.7, there exists $\lambda \in \mathbb{R}$ such that the spectral measure of $e^{i\lambda X_0^2} - \mathbb{E}[e^{i\lambda X_0^2}]$ is of the form $|\mathbb{E}[e^{i\lambda X_0^2}]|^2 e(\sigma^2)$, with $\sigma^2$ possessing an atom at 0. The spectral measure of $e^{i\lambda X_0^1} - \mathbb{E}[e^{i\lambda X_0^1}]$ is $|\mathbb{E}[e^{i\lambda X_0^1}]|^2 e(\sigma^1)$ for a measure $\sigma^1$. An easy computation shows that the spectral measure of $e^{i\lambda X_0} - \mathbb{E}[e^{i\lambda X_0}] = e^{i(\lambda X_0^1 + \lambda X_0^2)} - \mathbb{E}[e^{i(\lambda X_0^1 + \lambda X_0^2)}]$ is $|\mathbb{E}[e^{i\lambda X_0}]|^2 e(\sigma^1 + \sigma^2)$ but since $\sigma^1$ has an atom at 0, so has $|\mathbb{E}[e^{i\lambda X_0}]|^2 e(\sigma^1 + \sigma^2)$ and the process is not ergodic. Then, if $\mathbb{P}$ is ergodic, $\mathbb{P}_{Q_{ne}}$ is trivial and $\mathbb{P}$ writes $\mathbb{P}_{Q_B} * \mathbb{P}_{Q_m} * \mathbb{P}_{Q_{wm}}$ which implies that $\mathbb{P}$ is weakly mixing as a factor of the direct product of weakly mixing systems.   □

From Theorem 5.7, the hierarchy of "mixing" properties among ergodic IDp processes is explicit. Those process with a dissipative Lévy measure possess the strongest mixing behavior.

## 6. Generalized moving averages IDp processes.

DEFINITION 6.1.   A stationary process $\{X_n\}_{n \in \mathbb{Z}}$ is called *generalized moving average* if there exists an i.i.d. collection of processes $\{\{\xi_k^n\}_{k \in \mathbb{Z}}\}_{n \in \mathbb{Z}}$ such that, in distribution,

$$\{X_n\}_{n \in \mathbb{Z}} = \left\{ \sum_{k \in \mathbb{Z}} \xi_k^{n-k} \right\}_{n \in \mathbb{Z}}.$$

The process $\{\xi_k^0\}_{k \in \mathbb{Z}}$ is a *generator* of $\{X_n\}_{n \in \mathbb{Z}}$.

THEOREM 6.2.   *A stationary IDp process is generalized moving average with ID generator if and only if its Lévy measure is dissipative.*

Proof. Consider the distribution $\mathbb{P}$ of a generalized moving average IDp process of Lévy measure $Q$, the distribution $\mathbb{P}_g$ of an ID generator for it and $Q_g$ its Lévy measure. Since the process is the sum of the translates of independent process of distribution $\mathbb{P}_g$, we have

$$\mathbb{P} = \prod_{k \in \mathbb{Z}} \mathbb{P}_g \circ T^{-k}$$

(the product is the convolution) and thus,

$$Q = \sum_{k \in \mathbb{Z}} Q_g \circ T^{-k}.$$

We will show $Q$ is dissipative. Form the space $(\mathbb{Z} \times \mathbb{R}^{\mathbb{Z}}, \mathcal{P}(\mathbb{Z}) \otimes \mathcal{B}^{\otimes \mathbb{Z}}, m_c \otimes Q_g, \widetilde{T})$ where $m_c$ is the counting measure on $\mathbb{Z}$ and $\widetilde{T}$ is defined by $\widetilde{T}(n, \{x_i\}_{i \in \mathbb{Z}}) = (n+1, \{x_i\}_{i \in \mathbb{Z}})$, this system is clearly dissipative. Consider the map $\varphi$ from $\mathbb{Z} \times \mathbb{R}^{\mathbb{Z}}$ to $\mathbb{R}^{\mathbb{Z}}$ defined by $\varphi(n, \{x_i\}_{i \in \mathbb{Z}}) = T^n \{x_i\}_{i \in \mathbb{Z}}$. We have

$$\varphi \circ \widetilde{T}(n, \{x_i\}_{i \in \mathbb{Z}}) = \varphi(n+1, \{x_i\}_{i \in \mathbb{Z}})$$
$$= T^{n+1} \{x_i\}_{i \in \mathbb{Z}} = T(T^n \{x_i\}_{i \in \mathbb{Z}}) = T \circ \varphi(n, \{x_i\}_{i \in \mathbb{Z}})$$

and

$$(m_c \otimes Q_g) \circ \varphi^{-1} = \sum_{k \in \mathbb{Z}} Q_g \circ T^{-k}.$$

Thus, the map $\varphi$ is a factor map from $(\mathbb{Z} \times \mathbb{R}^{\mathbb{Z}}, \mathcal{P}(\mathbb{Z}) \otimes \mathcal{B}^{\otimes \mathbb{Z}}, m_c \otimes Q_g, \widetilde{T})$ to $(\mathbb{R}^{\mathbb{Z}}, \mathcal{B}^{\otimes \mathbb{Z}}, Q, T)$, this implies that $Q$ is dissipative.

Now assume that $Q$ is the dissipative Lévy measure of a stationary IDp process of distribution $\mathbb{P}$. From Lemma 2.6, there exists a wandering set $A$ such that $\mathbb{R}^{\mathbb{Z}} = \bigcup_{n \in \mathbb{Z}} T^{-n} A$ mod. $Q$. If we denote by $Q_g := Q_{|A}$ and $\mathbb{P}_g$ the distribution of the ID process with Lévy measure $Q_g$, then, since

$$Q = \sum_{k \in \mathbb{Z}} Q_g \circ T^{-k},$$

we obtain that

$$\mathbb{P} = \prod_{k \in \mathbb{Z}} \mathbb{P}_g \circ T^{-k}$$

and we can deduce that $\mathbb{P}_g$ is the distribution of an IDp process, generator for $\mathbb{P}$. □

## 7. Square integrable IDp processes.

Here we consider (with the exception of Proposition 7.2) square integrable IDp processes. To motivate this section, note that if $Q$ is a (shift)-stationary measure on $(\mathbb{R}^{\mathbb{Z}}, \mathcal{B}^{\otimes \mathbb{Z}})$ such that $\int_{\mathbb{R}^{\mathbb{Z}}} x_0^2 Q(dx) < +\infty$ satisfies $Q\{0\} = 0$, $Q$ can be considered as the Lévy measure of an IDp process which will prove to be square integrable. The family of Lévy measures of this type is hence quite large.

7.1. *Fundamental isometry.* We assume that the process is centered and we denote by $U$ (resp. $V$) the unitary operator associated to $T$ in $L^2(\mathbb{P})$ [resp. $L^2(Q)$] and $\mathfrak{C}_{X_0}(\mathbb{P})$ [resp. $\mathfrak{C}_{X_0}(Q)$] the cyclic subspace associated to $X_0$ in $L^2(\mathbb{P})$ [resp. $L^2(Q)$]. We establish the following result:

PROPOSITION 7.1. $\mathfrak{C}_{X_0}(\mathbb{P})$ *is unitary isometric to* $\mathfrak{C}_{X_0}(Q)$, *the unitary operators* $U$ *and* $V$ *being conjugate.*

PROOF. The property comes from the following identities:

$$\langle X_k, X_p \rangle_{L^2(\mathbb{P})} = \int_{\mathbb{R}^{\mathbb{Z}}} x_k x_p \mathbb{P}(dx)$$

$$= \int_{\mathbb{R}^2} uv \mathbb{P}_{(X_k, X_p)}(du, dv)$$

$$= \int_{\mathbb{R}^2} uv Q_{(X_k, X_p)}(du, dv)$$

$$= \int_{\mathbb{R}^{\mathbb{Z}}} x_k x_p Q(dx) = \langle X_k, X_p \rangle_{L^2(Q)}$$

(the equality between $\int_{\mathbb{R}^2} uv \mathbb{P}_{(X_k, X_p)}(du, dv)$ and $\int_{\mathbb{R}^2} uv Q_{(X_k, X_p)}(du, dv)$, where $Q_{(X_k, X_p)}$ denotes the Lévy measure of $(X_k, X_p)$, can be found in [19] page 163).

That is, if we denote by $\Phi$ the mapping that associates $X_k$ in $L^2(\mathbb{P})$ to $X_k$ in $L^2(Q)$ for all $k \in \mathbb{Z}$, then $\Phi$ can be extended linearly to an isometry between $\mathfrak{C}_{X_0}(\mathbb{P})$ and $\mathfrak{C}_{X_0}(Q)$. The fact that $\Phi U = V \Phi$ is obvious.

Thus, $X_0$ has the same spectral measure under $\mathbb{P}$ or under $Q$. □

7.2. *Ergodic and mixing criteria.* We recall the Gaussian case (see [2]), where ergodicity and mixing of the system is determined by the spectral measure of $X_0$:

- The system is ergodic if and only if $\sigma$ is continuous.
- The system is mixing if and only if $\sigma$ is Rajchman.

We then observe that, thanks to Proposition 7.1, such criteria no longer apply for square integrable IDp processes. Indeed, taking the distribution $Q$ of a centered square integrable mixing process, the IDp process with Lévy measure $Q$ is not ergodic by Theorem 5.7, but the spectral measure $\sigma$ of $X_0$ satisfies $\hat{\sigma}(k) \to 0$ as $|k|$ tends toward infinity. We must then assume some restrictions on the trajectories of the process to draw conclusions on ergodicity and mixing by only looking at the spectral measure of $X_0 - \mathbb{E}[X_0]$.

We start by a result where integrability suffices.

PROPOSITION 7.2. *Let $X$ be an IDp process of distribution $\mathbb{P}$ such that, up to a possible translation or a change of sign, $X_0$ is nonnegative. Then $\mathbb{P}$ is ergodic if and only if $\frac{1}{n}\sum_{k=1}^{n} X_k \to \mathbb{E}[X_0]$ $\mathbb{P}$-a.s. [or in $L^1(\mathbb{P})$] as $n$ tends to infinity.*

PROOF. We suppose that $X_0$ is nonnegative and that we have the representation (5.3) through the Laplace transform, $a \in \mathfrak{A} \cap \mathbb{R}_+^{\mathbb{Z}}$:

$$\mathbb{E}[\exp -\langle a, X \rangle] = \exp -\left[\int_{\mathbb{R}^{\mathbb{Z}}} 1 - e^{-\langle a, x \rangle} Q(dx)\right].$$

If one knows that $\mathbb{P}$ is ergodic, then the convergence holds thanks to the Birkhoff ergodic theorem.

Suppose now that $\frac{1}{n}\sum_{k=1}^{n} X_k \to \mathbb{E}[X_0]$ as $n$ tends to infinity $\mathbb{P}$-a.s. without ergodicity of $\mathbb{P}$. The decomposition of $\mathbb{P}$ is of the type $\mathbb{P}_e * \mathbb{P}_{Q_{ne}}$, where $\mathbb{P}_e$ is ergodic. Let $X^{ne}$ be of distribution $\mathbb{P}_{Q_{ne}}$ and $X^e$ be of distribution $\mathbb{P}_e$, assumed independent, such that $X^{ne} + X^e$ is of distribution $\mathbb{P}$.

The fact that $\frac{1}{n}\sum_{k=1}^{n}[(X^{ne} + X^e)_n] \to \mathbb{E}[X_0^{ne}] + \mathbb{E}[X_0^e]$ implies

$$\frac{1}{n}\sum_{k=1}^{n} X_n^{ne} \to \mathbb{E}[X_0^{ne}].$$

Hence, using

$$\mathbb{E}_{Q_{ne}}\left[\exp -\frac{1}{n}\sum_{k=1}^{n} X_k\right] = \exp -\left[\int_{\mathbb{R}^{\mathbb{Z}}} 1 - \exp\left[-\frac{1}{n}\sum_{k=1}^{n} x_k\right] Q_{ne}(dx)\right],$$

we note that the term of the left-hand side tends to $\exp -\mathbb{E}_{Q_{ne}}[X_0]$ by dominated convergence and, by continuity of the exponential, we then have

(7.1) $$\int_{\mathbb{R}^{\mathbb{Z}}} 1 - \exp\left[-\frac{1}{n}\sum_{k=1}^{n} x_k\right] Q_{ne}(dx) \to \mathbb{E}_{Q_{ne}}[X_0].$$

Under this representation, we also know, by (5.3), that

$$\mathbb{E}_{Q_{ne}}[X_0] = \int_{\mathbb{R}^{\mathbb{Z}}} x_0 Q_{ne}(dx).$$

Now consider the probability $\nu$ which is $T$-invariant and equivalent to $Q_{ne}$ and let $f := \frac{dQ_{ne}}{d\nu}$ ($f$ is $T$-invariant).

$fx_0$ is $\nu$-integrable and we can apply the Birkhoff ergodic theorem to deduce that

$$\frac{1}{n}\sum_{k=1}^{n} f \circ T^k x_k = f\left(\frac{1}{n}\sum_{k=1}^{n} x_k\right)$$

converges $\nu$-a.e. and in $L^1(\nu)$ to the conditional expectation of $fx_0$ with respect to the invariant $\sigma$-field which we denote by $\nu(fx_0|\mathcal{I})$. But, since $f$ is $T$-invariant and nonnegative, $\nu(fx_0|\mathcal{I}) = f\nu(x_0|\mathcal{I})$, that is, by dividing by $f$, $\frac{1}{n}\sum_{k=1}^n x_k$ converges $\nu$-a.e. to $\nu(x_0|\mathcal{I})$.

Since

$$\left(1 - \exp\left[-\frac{1}{n}\sum_{k=1}^n x_k\right]\right)f \leq f\left(\frac{1}{n}\sum_{k=1}^n x_k\right)$$

and by using the fact that $f(\frac{1}{n}\sum_{k=1}^n x_k)$ converges in $L^1(\nu)$, the sequence $(1 - \exp[-\frac{1}{n}\sum_{k=1}^n x_k])f$ is uniformly integrable and, since it tends $\nu$-a.e. to $(1 - \exp[-\nu(x_0|\mathcal{I})])$, we observe that

$$\int_{\mathbb{R}^{\mathbb{Z}}} 1 - \exp\left[-\frac{1}{n}\sum_{k=1}^n x_k\right] Q_{ne}(dx) = \int_{\mathbb{R}^{\mathbb{Z}}}\left(1 - \exp\left[-\frac{1}{n}\sum_{k=1}^n x_k\right]\right)f\nu(dx)$$

tends, as $n$ tends to infinity, to

$$\int_{\mathbb{R}^{\mathbb{Z}}} (1 - \exp[-\nu(x_0|\mathcal{I})])f\nu(dx).$$

But since $x_0 \geq 0$ and $Q_{ne}\{0\} = 0$ (and then $\nu\{0\} = 0$), we have $\nu(x_0|\mathcal{I}) > 0$ $\nu$-a.e., thus,

$$\int_{\mathbb{R}^{\mathbb{Z}}} (1 - \exp[-\nu(x_0|\mathcal{I})])f\nu(dx) < \int_{\mathbb{R}^{\mathbb{Z}}} \nu(x_0|\mathcal{I})f\nu(dx) = \int_{\mathbb{R}^{\mathbb{Z}}} x_0 f\nu(dx),$$

that is, the limit, as $n$ tends to infinity of $\int_{\mathbb{R}^{\mathbb{Z}}} 1 - \exp[-\frac{1}{n}\sum_{k=1}^n x_k]Q_{ne}(dx)$, is strictly less than $\int_{\mathbb{R}^{\mathbb{Z}}} x_0 Q_{ne}(dx)$. This contradicts (7.1), there is no term of the form $\mathbb{P}_{Q_{ne}}$ in the factorization of $\mathbb{P}$ and $\mathbb{P}$ is thus ergodic.  $\square$

We can now prove a proposition for square integrable processes:

PROPOSITION 7.3.  *Let $X$ be an IDp process of distribution $\mathbb{P}$ such that, up to a possible translation or a change of sign, $X_0$ is nonnegative. Let $\sigma$ be the spectral measure of $X_0 - \mathbb{E}[X_0]$.*

*$\mathbb{P}$ is ergodic if and only if $\sigma\{0\} = 0$.*

*$\mathbb{P}$ is mixing if and only if $\sigma$ is Rajchman.*

PROOF.  We know that $\sigma\{0\}$ equals the variance of $\mathbb{E}[X_0|\mathcal{I}]$. Moreover, the Birkhoff ergodic theorem tells us that $\frac{1}{n}\sum_{k=1}^n X_k \to \mathbb{E}[X_0|\mathcal{I}]$ $\mathbb{P}$-a.s. Thus, if $\sigma\{0\} = 0$, $\mathbb{E}[X_0|\mathcal{I}]$ is constant and equals $\mathbb{E}[X_0]$, so we can apply Proposition 7.2 to conclude. Now if $\sigma$ is Rajchman, by the isometry, $\sigma$ is also the spectral measure of $X_0$ under $Q$ and we get $\int_{\mathbb{R}^{\mathbb{Z}}} x_0 x_n Q(dx) \to 0$ as $n$ tends to infinity and we can apply the mixing criterion established by Rosiński and Żak in [16] (Corollary 3, page 282).

Both converse implications follow from Proposition 3.2.  $\square$

**8. $\alpha$-semi-stable and $\alpha$-stable processes.** We recall the definition of an $\alpha$-semi-stable (resp. $\alpha$-stable) distribution on $(\mathbb{R}, \mathcal{B})$. Denote by $D_b$ the application which associates $x \in \mathbb{R}$ to $bx \in \mathbb{R}$. Assume that $0 < \alpha < 2$.

DEFINITION 8.1. An *$\alpha$-semi-stable distribution of span $b$ $(b > 0)$* is an IDp distribution on $(\mathbb{R}, \mathcal{B})$ whose Lévy measure $\nu$ satisfies

$$\nu = b^{-\alpha} D_b^\star(\nu).$$

A distribution is said to be *$\alpha$-stable* if it is $\alpha$-semi-stable of span $b$ for all $b > 0$.

We will now discuss $\alpha$-semi-stable and $\alpha$-stable processes by introducing the application $S_b$ which associates $\{x_n\}_{n \in \mathbb{Z}} \in \mathbb{R}^\mathbb{Z}$ to $\{bx_n\}_{n \in \mathbb{Z}}$.

DEFINITION 8.2. A stationary process is said to be *$\alpha$-semi-stable of span $b$* if it is IDp and its Lévy measure $Q$ satisfies

$$(8.1) \qquad\qquad Q = b^{-\alpha} S_b^\star(Q).$$

A stationary process is said to be *$\alpha$-stable* if it is $\alpha$-semi-stable of span $b$ for all $b > 0$.

In particular, $S_b$ is nonsingular and commutes with the shift $T$. Remark that an $\alpha$-semi-stable distribution of span $b$ or an $\alpha$-semi-stable process of span $b$ is also $\alpha$-semi-stable of span $\frac{1}{b}$.

PROPOSITION 8.3. *The canonical factorization of Theorem 5.5 of an $\alpha$-semi-stable process of span $b$ is exclusively made of $\alpha$-semi-stable processes of span $b$.*

PROOF. It suffices to show that the $T$-invariant subsets of the partition given in the canonical decomposition of Proposition 2.11 are also $S_b$-invariant.

Consider $(\mathbb{R}^\mathbb{Z}, \mathcal{B}^{\otimes \mathbb{Z}}, Q, T)$, where $Q$ satisfies (8.1). Let $\mathfrak{P}$ be the part of type $\mathbf{II}_1$ of the system, then there exists a $T$-invariant function $f$ such that $\mathfrak{P} = \{f > 0\}$ and $\int_{\mathbb{R}^\mathbb{Z}} f \, dQ = 1$. Let $b > 0$. The function $f \circ S_b$ is $T$-invariant since $f \circ S_b \circ T = f \circ T \circ S_b = f \circ S_b$. Thus, from (8.1), $\int_{\mathbb{R}^\mathbb{Z}} f \circ S_b \, dQ = \int_{\mathbb{R}^\mathbb{Z}} f \, dS_b^\star(Q) = b^\alpha \int_{\mathbb{R}^\mathbb{Z}} f \, dQ = b^\alpha$, so the probability measure with density $b^{-\alpha} f \circ S_b$ with respect to $Q$ is $T$-invariant. Thus, $S_b^{-1} \mathfrak{P} = \{f \circ S_b > 0\} \subset \mathfrak{P}$. By the same arguments, $S_{1/b}^{-1} \mathfrak{P} \subset \mathfrak{P}$ and thus, $S_b^{-1}(S_{1/b}^{-1} \mathfrak{P}) \subset S_b^{-1} \mathfrak{P}$ and this shows $S_b^{-1} \mathfrak{P} = \mathfrak{P}$.

Now consider the $T$-invariant set $\mathcal{N}_+$ of Proposition 2.10. Let $A \subset \mathcal{N}_+$ be such that $0 < Q(A) < +\infty$. Then

$$Q((S_b^{-1}A) \cap T^{-k}(S_b^{-1}A)) = Q((S_b^{-1}A) \cap S_b^{-1}(T^{-k}A))$$
$$= Q(S_b^{-1}(A \cap T^{-k}A))$$
$$= b^{-\alpha}Q(A \cap T^{-k}A)$$

and thus, $\overline{\lim}_{k\to\infty} Q((S_b^{-1}A) \cap T^{-k}(S_b^{-1}A)) = \overline{\lim}_{k\to\infty} b^{-\alpha}Q(A \cap T^{-k}A) > 0$. Then, $S_b^{-1}A \subset \mathcal{N}_+$ so we have $S_b^{-1}\mathcal{N}_+ \subset \mathcal{N}_+$, and, by symmetric arguments, $S_b^{-1}\mathcal{N}_+ = \mathcal{N}_+$.

Consider $\mathfrak{D}$, the dissipative part of the system. From Lemma 2.6, there exists a wandering set $W$ such that $\mathfrak{D} = \bigcup_{n\in\mathbb{Z}} T^{-n}W$. Let $b > 0$ and consider the set $S_b^{-1}W$ [which is of nonzero $Q$-measure from (8.1)]. We have, if $n \neq m$, $T^{-n}(S_b^{-1}W) \cap T^{-m}(S_b^{-1}W) = \varnothing$; indeed, using the nonsingularity of $S_b$,

$$T^{-n}(S_b^{-1}W) \cap T^{-m}(S_b^{-1}W)$$
$$= S_b^{-1}(T^{-n}W) \cap S_b^{-1}(T^{-m}W)$$
$$= S_b^{-1}(T^{-n}W \cap T^{-m}W) = \varnothing.$$

Thus, $S_b^{-1}W$ is a wandering set, so $S_b^{-1}\mathfrak{D} \subset \mathfrak{D}$ since

$$S_b^{-1}\mathfrak{D} = S_b^{-1}\left(\bigcup_{n\in\mathbb{Z}} T^{-n}W\right)$$
$$= \bigcup_{n\in\mathbb{Z}} S_b^{-1}(T^{-n}W)$$
$$= \bigcup_{n\in\mathbb{Z}} T^{-n}(S_b^{-1}W)$$

and $\mathfrak{D}$ is, by definition, the union of all the wandering sets. We conclude $S_b^{-1}\mathfrak{D} = \mathfrak{D}$.

It is now easy to finish the proof by looking at the invariance of complements, intersections, and so on, and show the invariance of each set in the partition:

$$\mathfrak{D} \cup (\mathcal{C} \cap \mathcal{N}_0) \cup (\mathcal{N}_+ \cap \mathcal{N}) \cup \mathfrak{P}. \qquad \square$$

COROLLARY 8.4. *The canonical factorization of an $\alpha$-stable process is exclusively made of $\alpha$-stable processes.*

8.1. *$S\alpha S$-processes and factorizations.* The most frequently studied stationary $\alpha$-stable processes are the so-called $S\alpha S$-processes, where the distribution is preserved under the change of sign. In our framework, this means that the involution $S$ commutes with the shift $T$ and preserves the Lévy measure (note that $S$ commutes also with $S_b$). It is easy to see that the canonical factorization of an $S\alpha S$ process is only made of $S\alpha S$ processes.

We now show some connections existing between the decomposition of Theorem 5.5 and decompositions of an $S\alpha S$ process previously established respectively by Rosiński [13], Pipiras and Taqqu [12] and Samorodnitsky [18]. We first recall their results (we refer to these papers for precise definitions), the symbol "=" means "equality in distribution."

THEOREM 8.5.    *A stationary $S\alpha S$ process $X$ admits the unique following decomposition, where the sum is made of independent $S\alpha S$ processes:*
    *(Rosiński)*

$$X = X_r^1 + X_r^2 + X_r^3.$$

*$X_r^1$ is a mixed moving average process, $X_r^2$ is harmonizable, $X_r^3$ cannot be decomposed as the sum of a mixed moving average (or harmonizable) process and an independent $S\alpha S$ process:*
    *(Pipiras and Taqqu)*

$$X = X_{pt}^1 + X_{pt}^2 + X_{pt}^3 + X_{pt}^4.$$

*$X_{pt}^1$ is a mixed moving average process, $X_{pt}^2$ is harmonizable, $X_{pt}^3$ is associated to a cyclic flow without harmonizable component, $X_{pt}^4$ cannot be decomposed as the sum of a mixed moving average, or a harmonizable process or a process associated to a cyclic flow, and an independent $S\alpha S$ process.*
    *(Samorodnitsky)*

$$X = X_s^1 + X_s^2 + X_s^3.$$

*$X_s^1$ is a mixed moving average process, $X_s^2$ is associated to a conservative null flow, $X_s^3$ is associated to a positive flow.*

These authors study both discrete and continuous time in the same framework and, to avoid unnecessary different terminology, use "flow" to designate both an action of $\mathbb{R}$ and of $\mathbb{Z}$. There is a confusing terminology in the literature about null and positive flows (see the remark after Proposition 2.10) and here, Samorodnitsky uses the one found in Aaronson's book [1].

Here we recall that, in general, there can be an infinity of components in the decomposition, our criteria were mostly chosen with respect to the ergodic properties of the components. In that way, our decomposition is closer to Samorodnitsky's:

PROPOSITION 8.6.    *$X_s^1$ has a dissipative Lévy measure, $X_s^1 + X_s^2$ has a $\mathbf{II}_\infty$ Lévy measure and $X_s^3$ has a $\mathbf{II}_1$ Lévy measure.*

PROOF. Note that $X_s^1$, $X_r^1$ and $X_{pt}^1$ have the same distribution. Rosiński has shown it is a mixed moving average process which implies that it is a generalized moving average (with a $S\alpha S$ generator). By Theorem 6.2, its Lévy measure is dissipative. $X_s^1$ was proved to be mixing but, thanks to Theorem 5.7, it has indeed the Bernoulli property.

Samorodnitsky has shown that $X_s^1 + X_s^2$ is ergodic, thus, by Theorem 5.7, its Lévy measure is of type $\mathbf{II}_\infty$. The same author has also proved that there do not exist two independent $S\alpha S$ processes $Z_1$ and $Z_2$, one of them being ergodic and such that $X_s^3 = Z_1 + Z_2$ and this proves, in our framework, that the Lévy measure of $X_s^3$ is of type $\mathbf{II}_1$. Indeed, if we write the decomposition of Theorem 5.5 $X = X_B + X_m + X_{wm} + X_{ne}$, $X_{ne}$ has the same distribution as $X_s^3$ and $X_B + X_m + X_{wm}$ has the same distribution as $X_s^1 + X_s^2$.  $\square$

**9. Conclusion.** For the sake of simplicity, we have dealt with a single transformation, but many of the techniques used here can be applied more generally to the study of infinitely divisible objects whose Lévy measure is preserved by any kind of group actions, for example, the continuous time versions of our results are mostly straightforward, as are the multidimensional or the complex valued ones. The use of Poisson suspensions seems "natural" in some way.

For the interested reader, more ergodic oriented results can be found in the Ph.D. thesis of the author and, we hope, will be published soon.

**Acknowledgments.** First, the author would like to thank the anonymous referee for the very detailed remarks and advice that helped a lot to correct the first version of this paper. He is also very grateful to Nathalie Eisenbaum, Mariusz Lemańczyk, Krzysztof Frączek, Jean-Paul Thouvenot and David Mac Donald.

LABORATOIRE ANALYSE GÉOMÉTRIE
ET APPLICATIONS
UMR 7539
UNIVERSITÉ PARIS 13
99 AVENUE J. B. CLÉMENT
F-93430 VILLETANEUSE
FRANCE
E-MAIL: roy@math.univ-paris13.fr